%
%
%
%
%
%


\magnification\magstep1
\hsize = 148 true mm

\iffalse  

\input amssym
\def\bR{{\Bbb R}}               
\def\bP{{\Bbb P}}               
\def\bB{{\Bbb B}}               
\let\restr = \upharpoonright    
\let\stminus = \smallsetminus   

\else    

\def\bR{{\rm I\!R}}               
\def\bP{{\rm I\!P}}               
\def\bB{{\rm I\!B}}               
\let\restr = \mid          
\let\stminus = \setminus   %
\def\square{\hbox
         {\vrule height 6 pt depth 0 pt width 0,4 pt%
          \vbox to 6pt{\hrule height0,4pt width5,2pt
                       \kern 5,2pt
                       \hrule height0,4pt width5,2pt
                       \vss}%
          \vrule height 6 pt depth 0 pt width 0,4 pt}}
\def\beth{{\scriptstyle\rm beth}}  

\fi      

\def\CH {{\rm CH}}
\let\arrows = \rightarrow
\def\narrows {\not\rightarrow}
\def\cR{{\cal R}}               
\def\cP{{\cal P}}               
\def\cA{{\cal A}}               
\def\et{\mathop\&}              
\def\tu#1
  {\oalign{$#1$\crcr\hidewidth\vbox to.2ex{\hbox{\rm\char`\~}\vss}\hidewidth}}
\def\force{\mathrel\Vert\joinrel\mathrel-} 
\def\cl{\mathop{\rm cl}\nolimits} 
\def\cf{\mathop{\rm cf}\nolimits} 
\def\otp{\mathop{\rm otp}\nolimits} 
\def\Fr{\mathop{\rm Fr}\nolimits} 
\def\img{\mathchar"005C}          
\def\beginsection #1\par{\bigbreak\leftline{\bf#1}\nobreak\medskip}
\def\QED{\unskip\enspace\enspace$\square$\medbreak}


\centerline{On $\CH + 2^{\aleph_1}\arrows(\alpha)^2_2$
for $\alpha<\omega_2$}\bigskip
\centerline{by}\bigskip
\centerline{Saharon Shelah\footnote{$^{1)}$}
{Publication number 424. Partially supported by BSF.}}\bigskip

\beginsection 1. Introduction

We prove the consistency of
$$
\CH + 2^{\aleph_1}{\rm\ is\
arbitrarily\ large} + 2^{\aleph_1}\narrows(\omega_1\times\omega)^2_2
$$
(Theorem~1).  If fact, we can get
$2^{\aleph_1}\narrows[\omega_1\times\omega]^2_{\aleph_0}$,
see 1A.
In addition to this theorem, we give generalizations
to other cardinals (Theorems 2~and~3).  The $\omega_1\times\omega$ is
best possible as \CH\ implies
$$ 
\omega_3 \arrows (\omega\times n)^2_2. 
$$
We were motivated by the
question of Baumgartner [B1] on whether \CH\ implies
$\omega_3\arrows(\alpha)^2_2$ for $\alpha<\omega_2$ (if
$2^{\aleph_1}=\aleph_2$, it follows from the Erdos--Rado theorem). 
He proved the
consistency of positive answer with $\CH+2^{\aleph_1} > \aleph_3$, and
proved in~ZFC a related polarized partition relation (from \CH)
$$
{\aleph_3\choose\aleph_2}\arrows{\aleph_1\choose\aleph_1}^{1,1}_{\aleph_0}.
$$

{\bf Note: } The main proof here is that of Theorem~1. In that proof,
in the way things are set up, the main point is proving the
$\aleph_2$-c.c. The main idea in the proof is using $\bR$ (defined in
the proof). It turns out that we can use as elements of~$\cP$ (see the
proof) just pairs $(a,b)$. Not much would be changed if we used
$\langle\,(a_n,\alpha_n):n<\omega\,\rangle$, \ $a_n$~a good approximation
of the $n^{\rm th}$ part of the suspected monochromatic set of order type
$\omega_1\times\omega$. In 1A, 2~and~3 we deal with generalizations
and in Theorem~4 with complementary positive results.

\beginsection 2. The main result

{\bf Theorem 1. } Suppose
\item{(a)} \CH.
\item{(b)} $\lambda^{\aleph_1}=\lambda$.\par\noindent
Then there is an $\aleph_2$-c.c.\ \ $\aleph_1$-complete forcing notion
$\bP$ such that
\item{(i)} $|\bP| = \lambda$.
\item{(ii)} $\force_\bP\>\hbox{``} 2^{\aleph_1} = \lambda$, \
$\lambda\narrows(\omega_1\times\omega)^2_2\hbox{''}$.
\item{(iii)} $\force_\bP\CH$.
\item{(iv)} Forcing with $\bP$ preserves cofinalities and cardinalities.\par

{\bf Proof. } By Erdos and Hajnal [EH] there is an algebra $\bB$ with
$2^{\aleph_0}=\aleph_1$ \ $\omega$-place functions, closed under
composition (for simplicity only), such that
$$\displaylines{
\rlap{$\otimes$}\hfill
{\rm If\ } \alpha_n<\lambda{\rm\ for\ } n<\omega,{\rm\ then\ for\ some\ }k
\hfill\cr
\alpha_k\in\cl_\bB\{\,\alpha_l: k<l<\omega\,\}.
}$$
[$\otimes$ implies that for every large enough $k$, for every $m$, \
$\alpha_k\in\cl_\bB\{\,\alpha_k: m<l<\omega\,\}$.] \ Let
$$
\cR_\delta = \{\,b: b\subseteq\lambda,\ \otp(b)=\delta,\
\alpha\in b\Rightarrow b\subseteq\cl_\bB(b\stminus\alpha)\,\}.
$$
So by $\otimes$ we have
$$\displaylines{
\rlap{$\oplus$}\hfill
{\rm If\ \alpha\ is\ a\ limit\ ordinal,\ }
b\subseteq\lambda,\ \otp(b)=\alpha,\hfill\cr
\textstyle{\rm\ then\ for\ some\ }
\alpha\in b,\ \ b\stminus\alpha\in\bigcup_\delta\cR_\delta.
\cr
}$$
Let $\cR_{<\omega_1} = \bigcup_{\alpha<\omega_1}\cR_\alpha$. Let $\bP$ 
be the set of forcing conditions
$$(w,c,\cP)$$
where $w$ is a countable subset of $\lambda$, \ $c:[w]^2\arrows
\{{\rm red, green}\}=\{0,1\}$ (but we write $c(\alpha,\beta)$ instead
of $c(\{\alpha,\beta\})$),
and $\cP$ is a countable family of
pairs $(a,b)$ such that
\item{(i)} $a$, $b$ are subsets of $w$
\item{(ii)} $b\in\cR_{<\omega_1}$ and $a$ is a finite union of members
of~$\cR_{<\omega_1}$
\item{(iii)} $\sup(a)<\min(b)$
\item{(iv)} If $\sup(a) \le\gamma<\min(b)$, \ $\gamma\in w$, then
$c(\gamma,\cdot)$ divides $a$~or~$b$ into two infinite sets.\par

We use the notation
$$p = (w^p,c^p,\cP^p)$$
for $p\in\bP$. The ordering of the conditions is defined as follows:
$$
p\le q \iff w^p\subseteq w^q \et c^p\subseteq c^q \et \cP^p\subseteq\cP^q.
$$
Let
$$
\tu c = \bigcup\{\,c^p: p\in\tu G_\bP\,\}.
$$

{\bf Fact A. } $\bP$ is $\aleph_2$-complete.

{\bf Proof. } Trivial---take the union.\QED

{\bf Fact B. } For $\gamma<\lambda$, \ $\{\,q\in\bP:\gamma\in
w^q\,\}$ is open dense.

{\bf Proof. } Let $p\in\bP$. If $\gamma\in w^p$, we are done.
Otherwise we define $q$ as follows: $w^q=w^p\cup\{\gamma\}$, \
$\cP^q=\cP^p$, \ $c^q\restr w^p = c^p$ and $c^q(\gamma,\cdot)$ is
defined so that if $(a,b)\in\cP^q$, then $c^q(\gamma,\cdot)$ divides
$a$~and~$b$ into two infinite sets.\QED

{\bf Fact C. } $\force_\bP$ ``$2^{\aleph_1}\ge\lambda$ and
$\tu c:[\lambda]^2\arrows\{{\rm red, green}\}$''

{\bf Proof. } The second phrase follows from Fact B. For the first
phase, define $\tu \rho_\alpha\in{}^{\omega_1}2$, for $\alpha<\lambda$,
by: $\tu \rho_\alpha(i)=\tu c(0,\alpha+i)$. Easily
$$
\force_\bP \hbox{``}\tu \rho_\alpha\in{}^{\omega_1}2 {\rm\ and\ for\ }
\alpha<\beta<\lambda,\ \tu \rho_\alpha\neq\tu \rho_\beta\hbox{''};
{\rm\ so\ } \force_\bP\hbox{``}2^{\aleph_1}\ge\lambda\hbox{''}.
$$
\QED

{\bf Fact D. } $\bP$ satisfies the $\aleph_2$-c.c.

{\bf Proof. } Suppose $p_i\in\bP$ for $i<\aleph_2$. For each~$i$
choose a countable family $\cA^i$ of subsets of~$w^{p_i}$ such that
$\cA^i\subseteq\cR_{<\omega_1}$ and $(a,b)\in\cP^{p_i}$ implies
$b\in\cA^i$ and $a$ is a finite union of members of~$\cA^i$. For each
$\gamma\in c\in\cA^i$ choose a function $F^i_{\gamma,c}$ s.t.\
$F^i_{\gamma,c}(c\stminus(\gamma+1))=\gamma$. Let $v_i$ be the
closure of $w_i$ (in the order topology).

We may assume that $\langle\, v_i:i<\omega_2\,\rangle$ is a
$\Delta$-system (we have \CH) and that $\otp(v_i)$ is the same for all
$i<\omega_2$. W.l.o.g.\ for~$i<j$ the unique order-preserving function
$h_{i,j}$ from~$v_i$ onto~$v_j$ maps $p_i$ onto~$p_j$, \ $w^{p_i}\cap
w^{p_j} = w^{p_0}\cap w^{p_1}$ onto itself, and
$$F^i_{\gamma,c} = F^j_{h_{i,j}(\gamma),h_{i,j}\img c}$$
for $\gamma\in c\in\cA^i$ (remember: $\bB$ has $2^{\aleph_0}=\aleph_1$
functions only). Hence
$$\displaylines{
\rlap{$\otimes_1$}\hfill
h_{i,j} {\rm\ is\ the\ identity\ on\ } v_i\cap v_j {\rm\ for\ }i<j.
\hfill\cr
}$$
Clearly by the definition of $\cR_{<\omega_1}$ and the condition
on~$F^i_{\gamma,c}$: 
$$\displaylines{
\rlap{$\otimes_2$}\hfill
{\rm If\ }a\in\cA^i,\ i\ne j{\rm\ and\ }a\not\subseteq w^{p_i}\cap w^{p_j},
\hfill\cr
{\rm then\ }a\stminus(w^{p_i}\cap w^{p_j}){\rm\ is\ infinite.}\cr
}$$
We define $q$ as follows.
\item{} $w^q=w^{p_0}\cup w^{p_1}$.
\item{} $\cP^q=\cP^{p_0}\cup\cP^{p_1}$.
\item{} $c^p$ extends $c^{p_0}$ and $c^{p_1}$ in such a way that, for
$e\in\{0,1\}$,
\itemitem{$(*)$}  for every $\gamma\in w^{p_e}\stminus w^{p_{1-e}}$
    and every $a\in\cA^{1-e}$, \ $w^p(\gamma,\cdot)$ divides $a$ into
    two infinite parts, provided that
\itemitem{$(**)$} $a\stminus w^{p_e}$ is infinite.\par

This is easily done and $p_0\le q$, $p_1\le q$, {\it provided that\/}
$q\in\bP$. For this the problematic part is $c^q$ and, in particular,
part~(iv) of the definition of~$\bP$. So suppose $(a,b)\in\cP^q$,
e.g., $(a,b)\in\cP^{p_0}$. Suppose also $\gamma^*\in w^q$ so that
$\sup(a)\le\gamma^*<\sup(b)$. If $\gamma^*\in w^{p_0}$, there is no
problem, as $p_0\in\bP$. So let us assume $\gamma^*\in w^q\stminus
w^{p_0}=w^{p_1}\stminus w^{p_0}$. If $a\stminus w^{p_1}$ or
$b\stminus w^{p_1}$ is infinite, we are through in view of
condition~$(*)$ in the definition of $c^q$. Let us finally assume
$a\stminus w^{p_1}$ is finite. But $a\subseteq w^{p_0}$. Hence
$a\stminus(w^{p_0}\cap w^{p_1})$ is finite and $\otimes_2$~implies it
is empty, i.e.\ $a\subseteq w^{p_0}\cap w^{p_1}$. Similarly,
$b\subseteq w^{p_0}\cap w^{p_1}$. So $h_{0,1}\restr(a\cup b)$ is the
identity. But $(a,b)\in\cP^{p_0}$. But $h_{i,j}$ maps $p_i$ onto~$p_j$.
Hence $(a,b)\in\cP^{p_1}$. As
$p_1\in\bP$, we get the desired conclusion.\QED

{\bf Fact E. } $\force_\bP$ ``There is no $c$-monochromatic subset
of~$\lambda$ of order-type~$\omega_1\times\omega$.''

{\bf Proof. } Let $p$ force the existence of a counterexample. Let $G$
be $\bP$-generic over~$V$ with $p\in G$. In~$V[G]$ we can find
$A\subseteq\lambda$ of order-type~$\omega_1\times\omega$ such that
$c\restr[A]^2$~is constant. Let $A=\bigcup_{n<\omega}A_n$ where
$\otp(A_n)=\omega_1$ and $\sup(A_n)\le\min(A_{n+1})$. We can replace
$A_n$ by any $A'_n\subseteq A_n$ of the same cardinality. Hence we may
assume w.l.o.g.\ 
$$A_n\in\cR_{\omega_1}{\rm\qquad for\ }n<\omega.  \leqno(*)_1$$
Let $\delta_n=\sup(A_n)$ and
$$\eqalign{
\beta_n = \min\{\,\beta&:\delta\le\beta<\lambda,\ d(\beta,\cdot)
    {\rm\ does\ not}\cr
&\qquad{\rm divide\ }\bigcup\limits_{l\le n}A_l{\rm\ into\ two\
    infinite\ sets}\,\},\cr
}$$
where $d=\tu c^G$. Clearly $\beta_n\le\min(A_{n+1})$. Hence
$\beta_n<\beta_{n+1}$. Let $d_n\in\{0,1\}$ be such that
$d(\beta_n,\gamma)=d_n$ for all but finitely many
$\gamma\in\bigcup_{l\le n}A_l$. Let $u$ be an infinite subset
of~$\omega$ such that $\{\,\beta_n:n\in u\,\}\in\cR_\omega$. Let
$A_l=\{\,\alpha^l_i:i<\omega_1\,\}$ in increasing order. So $p$ forces
all this on suitable names
$$
\langle\,\tu\beta_n:n<\omega\,\rangle,\
\langle\,\tu\alpha^l_i:i<\omega_1\,\rangle,\
\langle\,\tu\delta_n:n<\omega\,\rangle.
$$

As $\bP$ is $\aleph_1$-complete, we can find $p_0\in\bP$ with
$p\le p_0$ so that $p_0$ forces $\tu\beta_l=\beta_l$ and
$\tu\delta_n=\delta_n$ for some $\beta_l$~and~$\delta_n$. We can
choose inductively conditions $p_k\in\bP$ such that $p_k\le p_{k+1}$
and there are $i_k<j_k$ and $\alpha^l_i$ (for $i<j_k$) with
$$\eqalignno{
p_{k+1} \force{}&\hbox{``}
    \alpha^l_{i_k}>\sup\{\,i:\tu\alpha^l_i\in w^{p_k}\,\},\cr
    &\alpha^l_i\in w^{p_{k+1}} {\rm\ for\ } i<j_k,\cr
    &\{\,\alpha^l_i:i<i_k\,\}\subseteq\cl_\bB\{\,\alpha^l_i:i_k<i<j_k\,\},\cr
    &\tu\alpha^l_i=\alpha^l_i{\rm\ for\ }i<j_k{\rm\ and}\cr
    &\gamma\in[\delta_m,\beta_m)\cap w^{p_k}{\rm\ implies\ }
    \tu c(\gamma,\cdot)\cr
    &\qquad{\rm divides\ }\{\,\alpha^l_i:i<j_k,\ l\le m\,\}{\rm\ into}\cr
    &\qquad\hbox{\rm two infinite sets.''}\cr
}$$
(remember our choice of~$\beta_m$). Let
$$\eqalignno{
l(*)&= \min(u)\cr
a&= \{\,\alpha^l_i:l\le l(*),\ i<\bigcup_k j_k\,\}\cr
b&= \{\,\beta_l:l\in u\,\}\cr
q&= (\bigcup_k w^{p_k},\bigcup_k c^{p_k},\bigcup_k\cP^{p_k}\cup\{(a,b)\}).\cr
}$$
Now $q\in\bP$. To see that $q$ satisfies condition~(iv) of the
definition of~$\bP$, let $\sup(a)\le\gamma<\min(b)$. Then
$\sup\{\,\alpha^{l(*)}_{i_k}:k<\omega\,\}\le\gamma<\beta_{l(*)}$. But
$\gamma\in w^p=\bigcup_k w^{p_k}$, so for some~$k$, \ $\gamma\in
w^{p_k}$. This implies
$$
\gamma\notin\left(\alpha^{l(*)}_{i_{k+1}},\delta_{l(*)}\right),
$$
whence $\gamma\ge\delta_{l(*)}$ and
$$
\{\,\alpha^l_i:l\le l(*),\ i<j_k\,\}\subseteq a,
$$
which implies the needed conclusion.

Also  $q\ge p_k\ge p$. But now, if $r\ge q$ forces a value to
$\alpha^{l(*)}_{\mathop\cup_k j_k}$; we get a contradiction.\QED

{\bf Remark 1A. } Note that the proof of Theorem~1 also gives the
consistency of $\lambda\narrows[\omega_1\times\omega]^2_{\aleph_0}$:
replace ``$c(\gamma,\cdot)$ divides a set~$x$ into two infinite
parts'' by ``$c(\gamma,\cdot)$ gets all values on a set~$x$.''

\beginsection 3. Generalizations to other cardinals

How much does the proof of Theorem~1 depend on~$\aleph_1$? Suppose we
replace $\aleph_0$ by~$\mu$.
    \medbreak
{\bf Theorem 2. } Assume $2^\mu=\mu^+<\lambda=\lambda^\mu$ and
$2\le\kappa\le\mu$. Then for
some $\mu^+$-complete $\mu^{++}$-c.c.\ forcing notion~$\bP$ of
cardinality $2^\mu$:
$$
\force_\bP 2^\mu = \lambda,\qquad
\lambda\narrows[\mu^+\times\mu]^2_\kappa.
$$

{\bf Proof. } Let $\bB$ and $\cR_\delta$ be defined as above (for
$\delta\le\mu^+$). Clearly
$$\hbox{%
If $a\subseteq\lambda$ has no last element, then for some $\alpha\in
a$, \ $a\stminus\alpha\in\bigcup_\delta\cR_\delta$.
}\leqno\oplus$$
Hence, if $\delta=\otp(a)$ is additively indecomposable, then
$a\stminus\alpha\in\cR_\delta$ for some $\alpha\in a$.

Let $\bP_\mu$ be the set of forcing conditions
$$(w,c,\cP)$$
where $w\subseteq\lambda$, \ $|w|\le\mu$, \ $c:[w]^2\arrows\{{\rm
red, green}\}$, and $\cP$ is a set of~$\le\mu$ pairs $(a,b)$ such
that
\item{(i)} $a$, $b$ are subsets of~$w$.
\item{(ii)} $b\in\cR_\mu$, and $a$ is a finite union of members
    of~$\bigcup_{\mu\le\delta<\mu^+}\cR_\delta$.
\item{(iii)} $\sup(a)<\min(b)$.
\item{(iv)} If $\sup(a)\le\gamma<\min(b)$, \ $\gamma\in w$, then the
function~$c(\gamma,\cdot)$ gets all values ($<\kappa$) on~$a$ or on~$b$.
\par

With the same proof as above we get
$$\eqalignno{
&\bP_\mu{\rm\ satisfies\ the\ \mu^{++}\hbox{\rm-}c.c.,}\cr
&\bP_\mu{\rm\ is\ \mu^+\hbox{\rm-}complete,}\cr
\noalign{\leftline{(so cardinal arithmetic is clear) and}}
&\force_{\bP_\mu}\ \lambda\narrows[\mu\times\mu]^2_\kappa.
}$$
$\square$\medbreak

What about replacing $\mu^+$ by an inaccessible $\theta$? We can
manage by demanding
$$\eqalign{
\{\,a\cap(\alpha,\beta)&:(a,b)\in\cP,\
    \bigcup_n\otp(a\cap(\alpha,\beta))\times n=\otp(a)\cr
    &\quad(\alpha,\beta){\rm\ maximal\ under\ these\ conditions}\,\}\cr
}$$
is free (meaning there are pairwise disjoint end segments) and by
taking care in defining the order. Hence the completeness drops to
$\theta$-strategical completeness. This is carried out in Theorem~3 below.
    \medbreak
{\bf Theorem 3. } Assume $\theta=\theta^{<\theta}>\aleph_0$ and
$\lambda=\lambda^{<\theta}$. Then for some $\theta^+$-c.c.\
$\theta$-strategically complete forcing~$\bP$, \ $|\bP|=\lambda$ and
$$
\force_\bP 2^\theta=\lambda,\ \lambda\narrows(\theta\times\theta)^2_2.
$$

{\bf Proof. } For $W$ a family of subsets of~$\lambda$, each with no
last element, let
$$\eqalign{
\Fr(W) = \{\, f &: f{\rm\ is\ a\ choice\ function\ on\ {\mit W}\ s.t.}\cr
&\quad\{\,a\stminus f(a): a\in W\,\}{\rm\ are\ pairwise\ disjoint}\,\}.\cr
}$$
If $\Fr(W)\ne\emptyset$, \ $W$ is called {\it free}.

Let $\bP_{<\theta}$ be the set of forcing conditions
$$(w,c,\cP,W)$$
where $w\subseteq\lambda$, \ $|w|<\theta$, \ $c:[w]^2\arrows\{{\rm
red, green}\}$, \ $W$~is a free family of~$<\theta$ subsets of~$w$,
each of which is in $\bigcup_{\delta<\theta}\cR_\delta$, and $\cP$ is
a set of~$<\theta$ pairs $(a,b)$ such that
\item{(i)} $a$, $b$ are subsets of~$w$.
\item{(ii)} $b\in\cR_\omega$.
\item{(iii)} $\sup(a)<\min(b)$ and for some
$\delta_0<\delta_1<\cdots<\delta_n$, \ $\delta_0<\min(a)$, \
$\sup(a)\le\delta_n$, \ $a\cap[\delta_l,\delta_{l+1})\in W$.
\item{(iv)} If $\sup(a)\le\gamma<\min(b)$, \ $\gamma\in w$, then
$c(\gamma,\cdot)$ divides $a$ or $b$ into two infinite sets.\par

We order $\bP_{<\theta}$ as follows:
$$\eqalignno{
p\le q{\rm\ \ if{}f\ \ }
&w^p\subseteq w^q,\ c^p\subseteq c^q,\ \cP^p\subseteq\cP^q,\
W^p\subseteq W^q\cr
&{\rm and\ every\ }f\in\Fr(W^p){\rm\ can\ be\ extended}\cr
&{\rm to\ a\ member\ of\ }\Fr(W^q).\cr
}$$
$\square$\medbreak


\beginsection 4. A provable partition relation

{\bf Claim 4. } Suppose $\theta>\aleph_0$, \ $n,r<\omega$ and
$\lambda=\lambda^{<\theta}$. Then
$$
(\lambda^+)^{r}\times n\arrows(\theta\times n,\theta\times r)^2_2.
$$

{\bf Proof. } We prove this by induction on~$r$. Clearly the claim
holds for $r=0, 1$. So w.l.o.g.\ we assume $r\ge2$. Let $c$ be a
2-place function from $(\lambda^+)^{r}\times n$ to $\{\rm red,
green\}$. Let $\chi=\beth_2(\lambda)^+$. Choose by induction on~$l$ a
model~$N_l$ such that
$$
N_l \prec (H(\chi),{\in},{<^*}),
$$
$|N_l|=\lambda$, \ $\lambda+1\subseteq N_l$, \ $N^{<\theta}_l\subseteq
N_l$, \ $c\in N_l$ and $N_l\in N_{l+1}$. Here $<^*$ is a well-ordering
of~$H(\chi)$. Let
$$
A_l =
\left[(\lambda^+)^{r}\times l,\ (\lambda^+)^{r}\times(l+1)\right),
$$
and let $\delta_l\in A_l\stminus N_l$ be such that $\delta_l\notin x$
whenever $x\in N_l$ is a subset of~$A_l$
and $\otp(x)<(\lambda^+)^{r-1}$. W.l.o.g.\ we
have $\delta_l\in N_{l+1}$. Now we shall show
$$\leqalignno{
&{\rm If}\ Y\in N_0,\ Y\subseteq A_m,\ |Y|=\lambda^+{\rm\ and\ }
\delta_m\in Y,& (*)\cr
&{\rm then\ we\ can\ find\ }\beta\in Y{\rm\ such\ that\
}c(\beta,\delta_m)={\rm red}.\cr
}$$

{\it Why $(*)$ suffices?\/} Assume $(*)$ holds. We can construct by
induction on~$i<\theta$ and for each~$i$ by induction on~$l<n$ an
ordinal~$\alpha_{i,l}$ s.t.\
\item{(a)} $\alpha_{i,l}\in A_l$ and
    $j<i\Rightarrow\alpha_{j,l}<\alpha_{i,l}$.
\item{(b)} $\alpha_{i,l}\in N_0$.
\item{(c)} $c(\alpha_{i,l},\delta_m)={\rm red}$\quad for $m<n$.
\item{(d)} $c(\alpha_{i,l},\alpha_{i_1,l_1})={\rm red}$\quad when $i_1<i$
or $i_1=i\et l_1<l$.\par

Accomplishing this suffices as $\alpha_{i,l}\in A_l$ and
$$
l<m \Rightarrow \sup A_l \le \min A_m.
$$
Arriving in the inductive process at~$(i,l)$, let
$$
Y=\{\,\beta\in A_l: c(\beta,\alpha_{j,m})={\rm red\quad if\ }
j<i,\ m<n,\ {\rm or\ }j=i,\ m<l\,\}.
$$
Now clearly $Y\subseteq A_\lambda$. Also $Y\in N_0$ as all parameters are
from~$N_0$, their number is~$<\theta$ and $N_0^{<\theta}\subseteq N_0$. Also
$\delta_l\in Y$ by the induction hypothesis (and $\delta_l\in A_l$).
So by $(*)$ we can find $\alpha_{i,l}$ as required.

{\bf Proof of $(*)$: } $Y\not\subseteq N_0$, because $\delta_m\in Y$
and $Y\in N_0$. As $|Y|=\lambda^+$, we have $\otp(Y)\ge\lambda^+$. But
$\lambda^+\arrows(\lambda^+,\theta)^2$, so there is $B\subseteq A_m$
s.t.\ $|B|=\lambda^+$ and $c\restr B\times B$ is constantly red or there
is $B\subseteq A_m$
s.t.\ $|B|=\theta$ and $c\restr B\times B$ is constantly green. In the
former case we get the conclusion of the claim. In the latter case we
may assume $B\in N_0$, hence $B\subseteq N_0$, and let $k\le n$ be
maximal s.t.\
$$
B'=\{\,\xi\in B:\bigwedge_{l<k}c(\delta_l,\xi)={\rm red}\,\}
$$
has cardinality $\theta$. If $k=n$, any member of~$B'$ is as required
in~$(*)$. So assume $k<n$. Now $B'\in N_k$, since $B\in N_0\prec N_k$
and $\{N_l,A_l\}\in N_k$ and $\delta_l\in N_k$ for~$l<k$. Also
$$
\{\,\xi\in B':c(\delta_k,\xi)={\rm red}\,\}
$$
is a subset of~$B'$ of cardinality~$<\theta$ by the choice of~$k$. So
for some $B''\in N_0$, \
$c\restr\{\delta_k\}\times(B'\stminus B'')$ is constantly green (e.g.,
as $B'\subseteq N_0$, and $N_0^{<\theta}\subseteq N_0$). Let
$$\eqalignno{
Z&=\{\,\delta\in A_k:c\restr\{\delta\}\times(B'\stminus B'')
    {\rm\ is\ constantly\ green}\,\}\cr
\noalign{\leftline{and}}
Z'&=\{\,\delta\in Z:(\forall\alpha\in B'\stminus B'')
    (\delta<\alpha\Leftrightarrow\delta_k<\alpha)\,\}.\cr
}$$
So $Z\subseteq A_k$, \ $Z\in N_k$, \ $\delta_k\notin N_k$ and therefore
$\otp(Z)=\otp(A_k)=(\lambda^+)^{r}$. Note that $k\ne l\Rightarrow
Z'=Z$ and $k=l\Rightarrow Z'=Z\stminus\sup(B'\stminus B'')$, so $Z'$
has the same properties. Now we apply the induction hypothesis: one of
the following holds (note that we can interchange the colours):
(a)~there is $Z''\subseteq Z'$, \ $\otp(Z'')=\theta\times n$, \
$c\restr Z''\times Z''$ is constantly red, wlog $Z''\in N_k$, or
(b)~there is $Z''\subseteq Z'$, \ $\otp(Z'')=\theta\times(r-1)$, \
$c\restr Z''\times Z''$ green and wlog $Z''\in N_k$.
If (a), we are done; if (b), \ $Z''\cup(B'\stminus B'')$ is as
required.\QED
    \medbreak
{\bf Remark 4A. } So $(\lambda^{+})^{n+1}\arrows(\theta\times n)^2$ for
$\lambda=\lambda^{<\theta}$, \ $\theta=\cf(\theta)>\aleph_0$ \ (e.g.,
$\lambda=2^{<\theta}$).
    \medbreak
{\bf Remark 4B. } Suppose $\lambda=\lambda^{<\theta}$, \
$\theta>\aleph_0$. If $c$ is a 2-colouring of
$(\lambda^{+r})^{s}\times n$ by~$k$~colours and every subset of
it of order type $(\lambda^{+(r-1)})^{s}\times n$ has a
monochromatic subset of order type~$\theta$ for each of the colours,
one of the colours being red, then by the last proof we get
\item{(a)} There is a monochromatic subset of order type $\theta\times
n$ and of colour red {\it or\/}
\item{(b)} There is a colour~$d$ and a set~$Z$ of order type
$(\lambda^{+r})^{s}$ and a set~$B$ of order type $\theta$ s.t.\
$B<Z$ or $Z<B$ and
$$
\{\,(\alpha,\beta):\alpha\in B,\ \beta\in Z{\rm\ or\
}\alpha\ne\beta\in B\,\}
$$
are all coloured with $d$.

So we can prove that for 2-colourings by $k$ colours~$c$
$$
(\lambda^{+r})^{s}\times n\arrows
(\theta\times n_1,\ldots,\theta\times n_k)^2
$$
when $r$, $s$, $n$ are sufficiently large (e.g.,
$n\ge\min\{\,n_l:l=1,\ldots,k,\ s\ge\sum_{l=1}^k n_l\,\}$ by
induction on~$\sum_{l=1}^k n_l$.

Note that if $c$ is a 2-colouring of~$\lambda^{+2k}$, then for some
$l<k$ and $A\subseteq\lambda^{+2k}$ of order type $\lambda^{+(2l+2)}$
we have
\item{$(*)$} If $A'\subseteq A$, $\otp(A')=\lambda^{+2l}$, and $d$ is
a colour which appears in~$A$, then there is $B\subseteq A'$ of order
type~$\theta$ s.t.\ $B$~is monochromatic of colour~$d$.\par

We can conclude $\lambda^{+2k}\arrows(\theta\times n)^2_k$.

%
%

\beginsection References

\noindent [B1]\quad J. Baumgartner, ??

\noindent [EH]\quad P. Erdos and A. Hajnal, ??

\bigskip{\obeylines
  Institute of Mathematics
  The Hebrew University
  Jerusalem
  Israel
\bigskip
  Department of Mathematics
  Rutgers University
  New Brunswick, NJ
  USA
}

\bye